\documentclass{amsart}
\usepackage{amssymb}
\newcommand\cP{\mathcal{P}}

\usepackage[all]{xy}
\newcommand\datver[1]{\def\datverp%
 {\par\boxed{\boxed{\text{#1; Run: \today}}}}}

\newcommand\boxb[1]{\square_b}

\newcommand\wh{\widehat}

\numberwithin{equation}{section}

\newcommand\paperbody%
        {}

\newtheorem{lemma}{Lemma}
\newtheorem{proposition}{Proposition}

\newtheorem{non-theorem}{Non-Theorem}

\theoremstyle{remark}

\newtheorem{remark}{Remark}

\newcommand\WF{\operatorname{WF}}

\newcommand\Tr{\operatorname{Tr}}

\newcommand\Fg{\mathfrak{g}}

\newcommand\cA{\mathcal{A}}

\newcommand\hE{\widehat{E}}
\newcommand\hbbE{\widehat{\bbE}}

\newcommand\bbE{\mathbb E}

\newcommand\bbQ{\mathbb Q}
\newcommand\bbR{\mathbb R}

\newcommand\bbZ{\mathbb Z}

\newcommand\CIc{{\mathcal{C}}^{\infty}_c}

\newcommand\CI{{\mathcal{C}}^{\infty}}

\newcommand\CmI{{\mathcal{C}}^{-\infty}}

\newcommand\Diag{\operatorname{Diag}}

\newcommand\cFNs{{}^{\Phi}\overline N\kern-1pt{}^*}

\newcommand\ind{\operatorname{ind}}
\newcommand\inda{\operatorname{ind_{\text{a}}}}

\newcommand\tr{\operatorname{tr}}

\newcommand\Hom{\operatorname{Hom}}
\newcommand\Id{\operatorname{Id}}

\newcommand\SU{\operatorname{SU}}
\newcommand\PU{\operatorname{PU}}

\newcommand\ci{${\mathcal{C}}^\infty$}

\newcommand\nul{\operatorname{null}}

\newcommand\supp{\operatorname{supp}}

\newcommand\Mand{\text{ and }}

\newcommand\Mfor{\text{ for }}

\newcommand\Mforsome{\text{ for some }}

\newcommand\Min{\text{ in }}

\datver{1.0A; Revised: 23-11-2006}
\begin{document}
\title[Equivariant and fractional index]
{Equivariant and fractional index of projective elliptic operators}

%lasteqno Fract-equiv-index@ 28

\author{V. Mathai}
\address{Department of Mathematics, University of Adelaide,
Adelaide 5005, Australia}
\email{vmathai@maths.adelaide.edu.au}
\author{R.B. Melrose}
\address{Department of Mathematics,
Massachusetts Institute of Technology,
Cambridge, Mass 02139}
\email{rbm@math.mit.edu}
\author{I.M. Singer}
\address{Department of Mathematics,
Massachusetts Institute of Technology,
Cambridge, Mass 02139}
\email{ims@math.mit.edu}

%\dedicatory{\sethouse\datverp}
\begin{abstract} In this note the fractional analytic index, for a
projective elliptic operator associated to an Azumaya bundle, of
\cite{MMS3} is related to the equivariant index of \cite{At74,Si73} for an
associated  transversally elliptic operator.
\end{abstract}
\maketitle

%\tableofcontents

\section*{Introduction}

Recall the setup in \cite{MMS3}. Let $\cA$ be an Azumaya bundle of rank $N$
over a compact oriented manifold $X$ and let $\cP$ denote the associated principal
$\PU(N)$-bundle of trivializations of $\cA.$ Let $\bbE=(E^+, E^-)$ denote a
pair of projective vector bundles associated to $\cA$ (or $\cP),$ which is
to say a projective $\bbZ_2$ superbundle. For each such pair, we defined in
\cite{MMS3}, projective pseudodifferential operators
$\Psi^\bullet_{\epsilon}(X,\bbE)$ with support in an
$\epsilon$-neighborhood of the diagonal in $X\times X.$ The principal
symbol $\sigma(D)$ of an elliptic projective pseudodifferential operator
$D,$ defines an element of the compactly-supported twisted K-theory
\begin{equation*}
[(\tau^*(\bbE),\sigma(D))] \in
K^0(T^*X, \tau^*\cA)
\label{Fract-equiv-index.18}\end{equation*}
where $\tau: T^*X \to X$ is the projection. The fractional analytic index of $D,$
which is defined using a parametrix, gives a homomorphism
\begin{equation}
\inda \colon K^0(T^*X, \tau^*\cA)\longrightarrow\bbQ.
\label{Fract-equiv-index.25}\end{equation}

On the other hand, the projective vector bundles $E^\pm$ can be realized as
vector bundles, $\hbbE=(\hE^+,\hE^-),$ in the ordinary sense over the total
space of $\cP$ with an action of $\wh G=\SU(N)$ which is equivariant with
respect to the action of $G=\PU(N)$ and in which the center, $\bbZ_N,$ acts
as the $N$th roots of unity. Following \cite{At74,Si73}, the $\wh
G$-equivariant pseudodifferential operators $\Psi^\bullet_{\wh G}(P,\hbbE)$
are defined for any equivariant bundles, as is the notion of 
transversal ellipticity. The principal symbol $\sigma(A)$ of a
transversally elliptic $\wh G$-equivariant pseudodifferential operator,
$A,$ fixes an element in equivariant K-theory
\begin{equation*}
[(\tau^*(\hbbE), \sigma(A))] \in K^0_{\wh G}(T^*_{\wh G} \cP)
\label{Fract-equiv-index.19}\end{equation*}
and all elements arise this way. The $\wh G$-equivariant index of $A$,
which is defined using a partial parametrix for $A,$ is a homomorphism, 
\begin{equation}
\ind_{\wh G} \colon K^0_{\wh G}(T^*_{\wh G} \cP)\longrightarrow \CmI(\wh G).
\label{Fract-equiv-index.21}\end{equation}
The restriction on the action of the center on the lift of a projective
bundle to $\cP,$ as opposed to a general $\wh G$-equivariant bundle for the
$\PU(N)$ action, gives a homomorphism,
\begin{equation}
\iota:K^0(T^*X,\tau^*\cA)\longrightarrow K^0_{\wh G}(T^*_{\wh G} \cP).
\label{Fract-equiv-index.22}\end{equation}

The diagonal action of $G$ on $\cP^2$ preserves the diagonal which therefore has a
basis of $G$-invariant neighborhoods. From the $\wh G$-equivariant
pseudodifferential operators, with support in a sufficiently small
neighborhood of the diagonal, there is a surjective pushforward
map, which is a homomorphism at the level of germs,
\begin{equation}
\pi_*
\colon
\Psi^\bullet_{\wh G,\epsilon}(\cP,\hbbE)\longrightarrow 
\Psi^\bullet_{\epsilon}(X,\bbE)
\label{Fract-equiv-index.20}\end{equation}
to the projective pseudodifferential operators. In fact this map preserves
products provided the supports of the factors are sufficiently close to the
diagonal. Moreover, pushforward sends transversally elliptic $\wh
G$-equivariant pseudodifferential operators to elliptic projective
pseudodifferential operators and covers the homomorphism
\eqref{Fract-equiv-index.22}
\begin{equation}
\iota[\sigma(\pi_*(A))]=[\sigma(A)] \Min K^0_{\wh G}(T^*_{\wh G} \cP)
\Mfor A\in \Psi^\bullet_{\wh G,\epsilon}(\cP,\bbE)\text{ elliptic.}
\label{Fract-equiv-index.24}\end{equation}
Proposition~\ref{13.10.2006.13} below relates these two pictures. Namely,
if $\phi\in\CI(\SU(N))$ has support sufficiently close to $e\in\SU(N)$ and
is equal to $1$ in a neighborhood of $e$ then the evaluation functional
${\rm ev}_\phi:\CmI(\wh G)\longrightarrow\bbR$ gives a commutative diagram
involving \eqref{Fract-equiv-index.25}, \eqref{Fract-equiv-index.21} and
\eqref{Fract-equiv-index.22}
\begin{equation}
\xymatrix{
K^0(T^*X, \tau^*\cA)\ar[d]^{ \inda }\ar[r]^{\iota}
&
K^0_{\wh G}(T^*_{\wh G} \cP)\ar[d]^{ \ind_{\wh G}}
\\
\bbQ
&
\CmI(\wh G)\ar[l]^{{\rm ev}_\phi}.
}
\end{equation}
Informally one can therefore say that the fractional analytic index, as
defined in \cite{MMS3}, is the coefficient of the delta function at the
identity in $\SU(N)$ of the equivariant index for transversally elliptic
operators on $\cP.$ Note that there may indeed be other terms in the
equivariant index with support at the identity, involving derivatives of
the delta function and there are terms supported at other points of $\bbZ_N.$ 
\footnote{We thank M. Karoubi for calling our attention to the omission of the
assumption of orientation in our reference to the Thom isomorphism in our
earlier paper \cite{MMS1}. It is unfortunate that we did not
explicitly reference his pioneering work on twisted K-theory and refer the
reader to his interesting new paper on the Arxiv, math/0701789 and the
references therein.}

\section{Transversally elliptic operators and the equivariant index}

As in \cite{At74,Si73}, let $X$ be a compact \ci\ manifold with a smooth
action of a Lie group, $G\ni g:X\longrightarrow X.$ In particular the
Lie algebra $\Fg$ of $G$ is realized as a Lie algebra of smooth vector fields
$L_a\in\CI(X;TX),$ $a\in\Fg,$ $[L_a,L_b]=L_{[a,b]}.$ Let $\Gamma \subset
T^*X$ denote the annihilator of this Lie algebra, so $\Gamma$ is also the
intersection over $G$ of the null spaces of the pull-back maps 
\begin{equation}
\Gamma \cap T^*_pX=\bigcap_{g\in G}\nul\left(g^*:T^*_pX\longrightarrow
T^*_{g^{-1}(p)}X\right).
\label{Fract-equiv-index.2}\end{equation}

Now, suppose that $\bbE=(E_+,E_-)$ is a smooth superbundle on $X$ which has a
smooth linear equivariant graded action of $G.$ Let $P\in\Psi^k(X;\bbE)$ be
a pseudodifferential operator which is invariant under the induced action
of $G$ on operators and which is transversally elliptic, that is its
characteristic variety does not meet $\Gamma:$
\begin{equation}
\Gamma\cap \Sigma (P)=\emptyset,\quad \Sigma (P)=\left\{\xi \in
T^*X\setminus0;\sigma(P)(\xi)\text{ is not invertible}\right\}.
\label{Fract-equiv-index.3}\end{equation}
Under these conditions (for compact $G)$ the equivariant index is defined
in \cite{At74,Si73} as a distribution on $G.$ In fact this can be done quite
directly. To do so, recall that for a function of compact support $\chi
\in\CIc(G),$ the action of the group induces a graded operator
\begin{equation}
T_{\chi}:\CI(X;\bbE)\longrightarrow\CI(X;\bbE),\quad
T_{\chi}u(x)=\int_{G}\chi(g)g^*udg.
\label{Fract-equiv-index.5}\end{equation}

\begin{proposition}\label{Fract-equiv-index.4} A transversally elliptic
pseudodifferential operator, $P,$ has a parametrix $Q,$ microlocally in a
neighborhood of $\Gamma$ and then for any $\chi\in\CIc(G),$ 
\begin{equation}
T_{\chi}\circ (PQ-\Id_-)\in\Psi^{-\infty}(X;\bbE_-)\Mand
T_{\chi}\circ (QP-\Id_+)\in\Psi^{-\infty}(X;\bbE_+)
\label{Fract-equiv-index.6}\end{equation}
are smoothing operators and 
\begin{equation}
\ind_G(P)(\chi)=\Tr(T_{\chi}(PQ-\Id_-))-\Tr(T_{\chi}(QP-\Id_+))
\label{Fract-equiv-index.7}\end{equation}
defines a distribution on $G$ which is independent of the choice of $Q.$
\end{proposition}

\begin{proof} The construction of parametrices is microlocal in any region
where the operator is elliptic, so $Q$ exists with the following constraint
on the operator wavefront set,
\begin{equation}
\left(\WF'(PQ-\Id_-)\cup\WF'(QP-\Id_+)\right)\cap\Gamma =\emptyset.
\label{Fract-equiv-index.9}\end{equation}
Thus $\WF'$ is the wavefront set of the Schwartz kernel of a
pseudodifferential operator, as a subset of the conormal bundle to the
diagonal which is then identified with the cotangent bundle of the manifold. 
Then $Q$ is unique microlocally in the sense that any other such
parametrix $Q'$ satisfies $\WF'(Q'-Q)\cap \Gamma =\emptyset.$ The definition
of $\Gamma$ means that for any pseudodifferential operator $A$ with
$\WF'(A)\cap\Gamma =\emptyset,$ $T_{\chi}A$ is smoothing and depends
continuously on $\chi.$ Thus \eqref{Fract-equiv-index.6} holds and by the
continuity of the dependence on $\chi$ defines a distribution on $G.$

To see the independence of the choice of parametrix, suppose that $Q_i,$
 $i=0,1$ are two choices. Then $Q_t=(1-t)Q_0+tQ_1$ is a homotopy of
 parametrices for $t\in[0,1]$ which defines a linear family of
 distributions with derivative 
\begin{multline*}
\frac{d}{dt}\left\{\Tr(T_{\chi}(PQ_t-\Id_-))-\Tr(T_{\chi}(Q_tP-\Id_+))\right\}\\
=\Tr(T_{\chi}(P(Q_1-Q_0)))-\Tr(T_{\chi}((Q_1-Q_0)P))\\
=\Tr(P(T_{\chi}(Q_1-Q_0))-\Tr((T_{\chi}(Q_1-Q_0))P)=0.
\label{Fract-equiv-index.11}\end{multline*}
Here we use the fact that $P$ is invariant under the action of $G$ and so
commutes with $T_{\chi}.$ The microlocal uniqueness of parametrices,
implies that $T_{\chi}(Q_1-Q_0)$ is a smoothing operator so the final line
follows from the vanishing of the trace on commutators where one factor is
pseudodifferential and the other is smoothing.
\end{proof}

\begin{proposition}\label{Fract-equiv-index.8} The distribution in
\eqref{Fract-equiv-index.7} reduces to the equivariant index of \cite{At74,Si73}.
\end{proposition}

\begin{proof} The Atiyah-Singer equivariant index for a transversally 
elliptic operator $P$ is equal to
$
\Tr(T_\chi(\Pi_0)) - \Tr(T_\chi(\Pi_1))
$
where $\Pi_j,$  $j=0,1$ are the orthogonal projections onto the nullspaces of
$P$ and $P^*$ respectively. The desired equality therefore involves only
an interchange of integrals, over $G$ and $X.$ Namely, if one chooses (by
averaging) a $G$-invariant parametrix $Q$ for $P,$ then the index in
\eqref{Fract-equiv-index.7},
$$
\ind_G(P) (\chi) = \Tr(T_\chi(PQ-I)) - \Tr(T_\chi(QP-I))
$$
is equal to
$
\Tr(T_\chi(\Pi_0)) - \Tr(T_\chi(\Pi_1)).
$
Let $K_\chi(x, y)$ denote the Schwartz kernel of the operator 
$T_{\chi},$ so $K_\chi(x, y) = \int_{G}  \delta_{gx}(y) \chi(g) dg.$ Thus
\begin{align*}
\Tr(T_\chi \circ \Pi_j) & = \int_{x\in X} \int_{y\in X} K_\chi(x, y) \tr 
(\Pi_j(y, x))  dy dx\\
& = \int_{x\in X} \int_{y\in X} \int_{G}  \delta_{gx}(y) \chi(g) dg \tr 
(\Pi_j(y, x))  dy dx\\
& = \int_{G}  \chi(g)  dg\int_{x\in X} \int_{y\in X}  \delta_{gx}(y) \tr 
(\Pi_j(y, x)) dy dx\\
& =  \int_{G}  \chi(g)  dg\int_{x\in X}  \tr(\Pi_j(gx, x)) dx\\
& =  \int_{G}  \chi(g) \,  {\rm char}(\Pi_j)(g)  dg
\end{align*}
which shows that
$$
\Tr(T_\chi(\Pi_0)) - \Tr(T_\chi(\Pi_1)) = \int_G \ind_G(P)(g) \chi(g) dg.
$$
\end{proof}

\begin{proposition}\label{Fract-equiv-index.12} Consider the subgroup of
$G$ defined by
\begin{equation}
G_f=\{g\in G;gx=x\Mforsome x\in X\}
\label{Fract-equiv-index.13}\end{equation}
then
\begin{equation}
\supp(\ind_{G}(P))\subset G_f.
\label{Fract-equiv-index.14}\end{equation}
\end{proposition}

\begin{proof} If $G_f=G$, then there is nothing to prove. 
Suppose that $G_f \ne G$. Then for $g\in G\setminus G_f$, the set $\{(gx,x);x\in X\}$ is disjoint
from the diagonal. It follows that if $\chi \in\CIc(G)$ has support
sufficiently close to $g$ and both $P$ and its parametrix $Q$ are chosen to have
Schwartz kernels with supports sufficiently close to the 
diagonal (which is always possible), then the supports of the Schwartz kernels of all
the terms in \eqref{Fract-equiv-index.6} are disjoint from the diagonal. It
follows that $\ind_{G}(P)(\chi)=0$ for such $\chi$ so
$g\notin\supp(\ind_{G}(P)).$
\end{proof}

\section{Fractional and equivariant index}

The finite central extension 
\begin{equation}
\bbZ_N\longrightarrow \SU(N)\longrightarrow \PU(N)
\label{13.10.2006.1}\end{equation}
is at the heart of the relation between the fractional and equivariant index.
From an Azumaya bundle over a compact, oriented smooth 
manifold $X$ we construct the principal $\PU(N)$-bundle
$\cP$ of trivializations. We will assume that the projective vector bundles in this 
section come equipped with a fixed hermitian structure.

\begin{lemma}\label{Fract-equiv-index.26} A projective vector bundle $E$ associated to
an Azumaya bundle $\cA$ over $X$ lifts to a vector bundle $\hE$ over $\cP$ with
an action of $\SU(N)$ which is equivariant with respect to the $\PU(N)$
action on $\cP$ and in which the center $\bbZ_N$ acts as the $N$th roots of
unity. 
\end{lemma}

Now, if $\bbE$ is a super projective vector bundle over $X,$ it lifts to a
super vector bundle $\hbbE$ over $\cP$ with $\SU(N)$ action. Consider the vector 
bundle $\hom(\hbbE)$ over $\cP$ of homomorphisms from $\hE^+$ to $\hE^-.$ Since
the action of $\SU(N)$ on $\hom(\hbbE)$ is by conjugation it descends to an
action of $\PU(N)$ and hence $\hom(\hbbE)$ descends to a vector bundle
$\hom(\bbE)$ on $X.$

The space $\Psi^m(\cP;\hbbE)$ of pseudodifferential operators over $\cP$
acting from sections of $\hE^+$ to $\hE^-$ may be identified with the
corresponding space of kernels on $\cP\times\cP$ which are distributional
sections of $\Hom(\hbbE)\otimes\Omega_R,$ the `big' homomorphism bundle
over $\cP^2$ with fiber at $(p,p')$ consisting of the homomorphisms from
$\hE^+_{p'}$ to $\hE^-_{p},$ tensored with the right density bundle and
with conormal singularities only at the diagonal. We are interested in the
$\SU(N)$-invariant part $\Psi^m_{\SU(N)}(\cP;\hbbE)$ corresponding to the
kernels which are invariant under the `diagonal' action of $\SU(N).$

\begin{proposition}\label{13.10.2006.2} If $\Omega\subset\cP^2$ is a
sufficiently small neighborhood of $\Diag\subset\cP^2$ invariant under
the diagonal $\PU(N)$-action there is a well-defined push-forward map
into the projective pseudodifferential operators
\begin{equation}
\left\{P\in\Psi^{m}_{\SU(N)}(\cP;\hbbE);\,\supp(P)\subset\Omega\right\}\ni A
\longrightarrow \pi_*(A)\in\Psi^m_\epsilon(X;\bbE)
\label{13.10.2006.3}\end{equation}
which preserves composition of elements with support in $\Omega'$ such
that $\Omega '\circ\Omega '\subset\Omega.$ 
\end{proposition}

\begin{proof} The push-forward map extends the averaging map in the
$\PU(N)$-invariant case in which the action of $\bbZ_N$ is trivial. Then
\begin{equation}
\pi^*(\pi_*(A)\phi) =A\pi^*\phi
\label{13.10.2006.4}\end{equation}
defines $\pi_*(A)$ unambiguously, since $\pi^*\phi$ is a $\PU(N)$-invariant
section and hence so is $A\pi^*\phi,$ so determines a unique section of
the quotient bundle. It is also immediate in this case that 
\begin{equation}
\pi_*(AB)=\pi_*(A)\pi_*(B)
\label{13.10.2006.6}\end{equation}
by the assumed $\PU(N)$-invariance of the operators. The definition
\eqref{13.10.2006.4} leads to a formula for the Schwartz kernel of
$\pi_*(A).$ Namely, writing $A$ for the Schwartz kernel of $A$ on $\cP^2,$  
\begin{equation}
\pi_*A(x,x')=\int_{\pi^{-1}(x)\times\pi^{-1}(x')} A(p,p').
\label{13.10.2006.5}\end{equation}
Since the projection map is a fibration, to make sense of this formal
integral we only need to use the fact that the bundle, of which the
integrand is a section, is naturally identified with the pull-back of a
bundle over the base tensored with the density bundle over the domain. The
composition formula \eqref{13.10.2006.6} then reduces to Fubini's theorem,
using the invariance of the kernels under the diagonal $\PU(N)$ action.

In the projective case we instead start from the formula
\eqref{13.10.2006.5}. As shown in \cite{MMS3}, the vector bundle $\hom(\bbE)$ over
$X$ lifted to the diagonal in $X^2$ extends to a small neighborhood
$\Omega$ of the diagonal as a vector bundle $\Hom(\bbE)$ with composition
property. In terms of the vector bundle $\Hom(\hbbE)$ over $\cP^2$ this can be
seen from the fact that each point in $(x,x')\in\Omega$ is covered by a
set of the form
\begin{equation}
\{(gp,g'p');g,\ g'\in\SU(N),\ g'g^{-1}\in B\}
\label{Fract-equiv-index.28}\end{equation}
for $B\subset\SU(N)$ some small neighborhood of the identity.
Namely if $p$ is any lift of $x$ then there is a lift $p'\in\cP$ of $x'$
which is close to $p$ and all such lifts are of the form
\eqref{Fract-equiv-index.28}. The diagonal action on $\Hom(\hbbE)$ descends
to a $\PU(N)$ action and it follows that $\Hom(\hbbE)$ may be naturally
identified over the set \eqref{Fract-equiv-index.28} with the fiber of
$\Hom(\bbE).$ Hence $\Hom(\hbbE)$ may be identified over a neighborhood of
the diagonal in $\cP^2$ with the pull-back of $\Hom(\bbE)$ and this
identification is consistent with the composition property.

Thus over the fiber of the push-forward integral \eqref{13.10.2006.5} the
integrand is identified with a distributional section of the bundle lifted
from the base. The properties of the push-forward, that it maps the
kernels of pseudodifferential operators to pseudodifferential operators and
respects products, then follow from localization, since this reduces the problem
to the usual case discussed initially.
\end{proof}

\begin{lemma}\label{13.10.2006.8} If $\chi\in\CI(\SU(N))$ is equal to $1$ in a
neighborhood of $e\in\SU(N)$ and $\Omega$ is a sufficiently small neighborhood
of the diagonal in $\cP^2,$ depending on $\chi,$ then under the
push-forward map of Proposition~\ref{13.10.2006.2}
\begin{equation}
\Tr(\pi_*(A))=\Tr(T_{\chi}A),\ A\in\Psi^{-\infty}_{\SU(N)}(\cP;\hbbE),\
\supp(A)\subset\Omega.
\label{13.10.2006.9}\end{equation}
\end{lemma}

\begin{proof} In a local trivialization of $\cP$ the kernel of $T_{\chi}A$
is of the form
\begin{equation}
\int_{\SU(N)}\chi(h^{-1}g)A(x,h,x',g')dh
\label{13.10.2006.10}\end{equation}
so the trace is 
\begin{equation*}
 \int_{\SU(N)\times\SU(N)}\chi(h^{-1}g)A(x,h,x,g)dhdgdx=
\int_{\SU(N)}\chi(h)A(x,h,x,e)dh
\label{13.10.2006.11}\end{equation*}
using the invariance of $A.$ Since the support of $A$ is close to the
diagonal and $\chi=1$ close to the identity, $\chi=1$ on the support
this reduces to 
\begin{equation*}
\int_{\SU(N)}A(x,h,x,e)dh=\Tr(\pi_*(A))
\label{13.10.2006.12}\end{equation*}
again using the $\SU(N)$-invariance of $A.$ 
\end{proof}

Now, suppose $A\in\Psi^m_{\SU(N)}(\cP,\bbE)$ is transversally
elliptic. Then the $\SU(N)$-equivariant index is the distribution 
\begin{equation}
\ind_{\SU(N)}(A)(\chi)=\Tr\left(T_{\chi}(AB-\Id_-))-\Tr(T_{\chi}(BA-\Id_+)\right)
\label{13.10.2006.7}\end{equation}
where $B \in\Psi^{-m}_{\SU(N)}(\cP,\bbE)$ is a parametrix for $A$ and $\chi
\in \CI(\SU(N)).$  We may choose $A$ and $B$ to have (kernels with)
supports arbitrarily close to the diagonal but maintaining the
$\SU(N)$-invariance.

\begin{proposition}\label{13.10.2006.13} If $\phi\in\CI(\SU(N))$ has support
sufficiently close to $e\in\SU(N)$ and is equal to $1$ in a neighborhood of
$e$ then  
\begin{equation}
\ind_{\SU(N)}(A)(\phi)=\inda({\pi_*(A)})
\label{13.10.2006.14}\end{equation}
for any transversally elliptic $A\in\Psi^m_{\SU(N)}(\cP;\hbbE)$ with support
  sufficiently close to the diagonal.
\end{proposition}

\begin{proof} Using Lemma~\ref{13.10.2006.8},
\begin{equation}
\Tr\left({T_\phi}(AB-\Id_-)\right)=\Tr(\pi_*(A)\pi_*(B)-\pi_*(\Id_+)),
\label{13.10.2006.15}\end{equation}
and similarly for the second term. Since $\pi_*(\Id)=\Id,$
$\pi_*(A)\pi_*(B)-\Id$ is a smoothing operator. In particular $\pi_*(B)$ is
a parametrix for $\pi_*(A)$ and \eqref{13.10.2006.14} follows.
\end{proof}

\begin{remark} Every compact, oriented, Riemannian manifold $X$ of
dimension $2n$, has projective vector bundles of half spinors, which are
realized as $\SU(N)$-equivariant vector bundles $\widehat{\mathbb S} =
(\widehat S^+, \widehat S^-),$ $N=2^n$, over the principal $\PU(N)$-bundle
$\cP$ over $X$ that is associated to the oriented orthonormal frame bundle
of $X$, cf\@. \S3 in \cite{MMS3}. Explicitly, $\widehat{\mathbb S}$ is the
$\bbZ_2$-graded $\SU(N)$-equivariant vector bundle of spinors associated to
the conormal bundle to the fibers, $T^*_{\SU(N)}\cP.$ On $\cP$ there is a
transversally elliptic, $\SU(N)$-equivariant Dirac operator $\eth^+,$
defined as follows. The Levi-Civita connection on $X$ determines in an
obvious way, partial spin connections $\nabla^\pm$ on $\widehat S^\pm$.
That is, $\nabla^+ : \CI(\cP, \widehat S^+) \to \CI(\cP, T^*_{\SU(N)}\cP
\otimes \widehat S^+)$.  If $C : \CI(\cP, T^*_{\SU(N)}\cP \otimes \widehat
S^+) \to \CI(\cP, \widehat S^-) $ denotes contraction given by Clifford
multiplication, then $\eth^+ : \CI(\cP, \widehat S^+) \to \CI(\cP, \widehat
S^-)$ is defined as the composition, $C\circ \nabla^+.$

Then $\pi_*(\eth^+)$ is just the projective Dirac operator
of \cite{MMS3}, and Proposition~\ref{13.10.2006.13} relates the indices in
these two senses.
\end{remark}

\begin{remark} Once the pushforward map $\pi_* \colon
\Psi^\bullet_{\wh G}(\cP,\hbbE)\longrightarrow
\Psi^\bullet_{\epsilon}(X,\bbE)$ is defined, Proposition
\ref{13.10.2006.13} can also be deduced from the index theorem of
\cite{MMS3} and the explicit topological expression for the equivariant
transversal index as in \cite{BV, BV2} simply by comparing the formul\ae.
\end{remark}

\providecommand{\href}[2]{#2}

\end{document}